\documentclass[9pt,shortpaper,twoside,web]{ieeecolor}
\usepackage{generic}
\usepackage{cite}
\usepackage{amsmath,amssymb,amsfonts}
\usepackage{algorithmic}
\usepackage[ruled]{algorithm}
\usepackage{graphicx}
\usepackage{pgfplots}
\usetikzlibrary{positioning}
\pgfplotsset{compat=1.14} 
\tikzstyle{vertex} = [fill,shape=circle,node distance=30pt]
\tikzstyle{edge} = [fill,opacity=.6,fill opacity=.5,line cap=round, line join=round, line width=10pt]
\tikzstyle{elabel} =  [fill,shape=circle,node distance=30pt,fill opacity=.9]
\definecolor{mygray}{gray}{0.95}
\pgfdeclarelayer{background}
\pgfsetlayers{background,main}
\usepackage[most]{tcolorbox}
\definecolor{mypurple}{rgb}{0.59, 0.44, 0.84}
\usepackage[mathscr]{euscript}

\usepackage{textcomp}

\begin{document}
\title{Explicit Solutions and Stability Properties of Homogeneous Polynomial Dynamical Systems}
\author{Can Chen
\thanks{This work is supported by Mathematics Department Graduate Fellowship from the University of Michigan. }
\thanks{Can Chen is with the Department of Mathematics, University of Michigan, Ann Arbor, MI 48109 USA and with the Channing Division of Network Medicine, Brigham and Women's Hospital and Harvard Medical School, Boston, MA 02115 USA (canc@umich.edu).}}

\maketitle

\begin{abstract}
In this paper, we provide a system-theoretic treatment of certain continuous-time homogeneous polynomial dynamical systems (HPDS) via tensor algebra. In particular, if a system of homogeneous polynomial differential equations can be represented by an orthogonally decomposable (odeco) tensor, we can construct its explicit solution by exploiting tensor Z-eigenvalues and Z-eigenvectors. We refer to such HPDS as odeco HPDS. By utilizing the form of the explicit solution, we are able to discuss the stability properties of an odeco HPDS.  We illustrate that the Z-eigenvalues of the corresponding dynamic tensor can be used to establish necessary and sufficient stability conditions, similar to these from linear systems theory. In addition, we are able to obtain the complete solution to an odeco HPDS with constant control. Finally, we establish results which enable one to determine if a general HPDS can be transformed to or approximated by an odeco HPDS, where the previous results can be applied. We demonstrate our framework with simulated and real-world examples. 
\end{abstract}


\section{Introduction}

Many real-world models such as those arising in biology, chemistry, and engineering can be captured by homogeneous polynomial dynamical systems (HPDS) \cite{chen2021controllability, bairey2016high,chellaboina2009modeling,craciun2006understanding,donnell2013local,smit2019walking,zhao2019optimal}. For example, Chen \textit{et al.} \cite{chen2021controllability} utilized HPDS to model  higher-order interactions in mouse neuronal networks, which offers a unique insight in understanding the functionality of  mouse brains. In addition, Bairey \textit{et al.} \cite{bairey2016high} used a polynomial dynamics model to simulate communities with interactions of different orders in ecological systems and analyzed the impacts of higher-order interactions on critical community size. Under certain conditions, the polynomial dynamics model can be simplified to HPDS. Investigating the stability properties of such higher-order networks is extremely significant in order to maintain and control the networks. Nevertheless, the stability of HPDS is one of the most challenging problems in systems theory due to its nature of nonlinearity. The most common strategy is still the local stability analysis when dealing with HPDS.

Numerous efforts have been made in exploring the stability properties of HPDS \cite{samardzija1983stability, ahmadi2019algebraic,ahmadi2013stability,ji2013constructing,she2013discovering,ahmadi2017sum}. In 1983, Samardzija \cite{samardzija1983stability} established a necessary and sufficient condition for asymptotic stability in two-dimensional HPDS by formulating a generalized characteristic value problem. Moreover, Ali and Khadir \cite{ahmadi2019algebraic} discovered that the existence of a rational Lyapunov function (i.e., the ratio of two polynomials) is necessary and sufficient for asymptotic stability of a HPDS. More importantly, they proved that the Lyapunov inequalities on both the rational function and its derivative have sum of squares certificates, so a Lyapunov function can always be found by semidefinite programming. The semidefinite programming problem depends on the two degree parameters, which gives rise to a hierarchy of semidefinite programs where one has to try all possible combinations of the parameters in order to obtain a Lyapunov function. Recently, tensor algebra has been applied to model and simulate nonlinear dynamics \cite{chen2021stability,chen2021controllability, gelss2017tensor,kruppa2017comparison,kruppa2018feedback,doi:10.1137/1.9781611975758.18}. Notably, every HPDS can be represented by a dynamic tensor, analogous to dynamic matrices of linear systems. Thus, tensor algebra such as tensor eigenvalues and tensor decompositions could be promising in determining the explicit solution and stability properties of a HPDS.

There are many definitions of tensor decompositions \cite{chen2021multilinear, Kolda06multilinearoperators, doi:10.1137/07070111X,doi:10.1137/140989340}. Of particular interest in this paper is tensor orthogonal decomposition. Generalized from matrix eigenvalue decomposition, tensor orthogonal decomposition aims to decompose a supersymmetric tensor (i.e., invariant under any permutations of the indices) into a sum of rank-one tensors in the form of outer products of vectors which form an orthonormal basis, see Fig. 1 B. Each rank-one tensor is also multiplied with a real coefficient. It has been shown that the coefficients and  the vectors in the orthogonal decomposition are the Z-eigenvalues and the Z-eigenvectors of the tensor, respectively \cite{doi:10.1137/140989340}. If a supersymmetric tensor has an orthogonal decomposition, it is called orthogonally decomposable (odeco). Odeco tensors enjoy the nice orthonormal property, which can be applied to various tensor-based applications. For instance, Anandkumar \textit{et al.} \cite{anandkumar2014tensor} exploited odeco tensors to estimate parameters in the method of moments from statistics. Thus, we are intrigued  by exploring the system properties of HPDS that can be represented by odeco tensors. We refer to such HPDS as odeco HPDS. Although the class of odeco HPDS is restricted, a certain amount of population dynamics with pairwise/higher-order interactions such as those arising in neuronal networks, chemical reaction networks, and ecological networks can be modeled by odeco HPDS. Most crucially, the results of odeco HPDS could be a powerful foundation for extension to general HPDS, e.g., transforming/approximating the dynamic tensor of a general HPDS to/by an odeco tensor.

In fact, Chen \cite{chen2021stability} investigated the explicit solutions and the stability properties of discrete-time odeco HPDS (also called multilinear dynamical systems in \cite{chen2021stability}). The author showed that the Z-eigenvalues from the orthogonal decomposition of a dynamic tensor play a significant role in the stability analysis, offering necessary and sufficient conditions. In this paper, we will focus on continuous-time HPDS. The key contributions of the paper are:
\begin{enumerate}
    \item We investigate the explicit solution of an odeco HPDS. We derive an explicit solution formula by using the Z-eigenvalues and the Z-eigenvectors from the orthogonal decomposition of the corresponding dynamic tensor.  
    \item According to the formula of the explicit solution, we are able to discuss the stability properties of an odeco HPDS. We show that the Z-eigenvalues from the orthogonal decomposition can offer necessary and sufficient stability conditions. Furthermore, we apply an upper bound of the largest Z-eigenvalue to determine the asymptotic stability efficiently.
    \item We explore the complete solution of an odeco HPDS with constant control. We find that the complete solution can be solved implicitly by exploiting Gauss hypergeometric functions. 
    \item We establish results which enable one to determine if a general HPDS can be transformed to or approximated by an odeco HPDS, in which all the previous results can be applied.  
\end{enumerate}

The paper is organized into seven sections. In Section \ref{sec:2}, we review tensor preliminaries including tensor vector multiplications, tensor eigenvalues, and orthogonal decomposition. We derive an explicit solution formula for an odeco HPDS and discuss the stability properties of the HPDS based on the form of the explicit solution in Section \ref{sec:3}. In Section \ref{sec:4}, we explore the complete solution of an odeco HPDS with constant control. We provide criteria to determine if a general HPDS can be transformed to or approximated by an odeco HPDS with detailed algorithmic procedures in Section \ref{sec:5}. Simulated and real-world examples are presented in Section \ref{sec:6}. Finally, we conclude in Section \ref{sec:7} with future research directions.

\section{Tensor Preliminaries}\label{sec:2}
A tensor is a multidimensional array \cite{chen2021multilinear, doi:10.1137/S0895479896305696,doi:10.1137/07070111X,9119161}. The order of a tensor is the number of its dimensions,  and each dimension is called a mode. A $k$th-order tensor is usually denoted by $\textsf{T}\in \mathbb{R}^{n_1\times n_2\times  \dots \times n_k}$.  It is therefore reasonable to consider scalars $x\in\mathbb{R}$ as zero-order tensors, vectors $\textbf{v}\in\mathbb{R}^{n}$ as first-order tensors, and matrices $\textbf{M}\in\mathbb{R}^{m\times n}$ as second-order tensors. A tensor is called cubical if every mode is the same size, i.e., $\textsf{T}\in \mathbb{R}^{n\times n\times  \dots \times n}$. A cubical tensor $\textsf{T}$ is called supersymmetric if $\textsf{T}_{j_1j_2\dots j_k}$ is invariant under any permutation of the indices. For instance, a third-order tensor $\textsf{T}\in\mathbb{R}^{n\times n\times n}$ is supersymmetric if 
\begin{equation*}
    \textsf{T}_{j_1j_2j_3}=\textsf{T}_{j_1j_3j_2}=\textsf{T}_{j_2j_1j_3}=\textsf{T}_{j_2j_3j_1}=\textsf{T}_{j_3j_1j_2}=\textsf{T}_{j_3j_2j_1}
\end{equation*}
for all $j_1,j_2,j_3=1,2,\dots,n$.

\subsection{Tensor Vector Multiplication}
The tensor vector multiplication $\textsf{T} \times_{p} \textbf{v}$ along mode $p$ for a vector $\textbf{v}\in  \mathbb{R}^{n_p}$ is defined as
\begin{equation}
(\textsf{T} \times_{p} \textbf{v})_{j_1j_2\dots j_{p-1}j_{p+1}\dots j_k}=\sum_{j_p=1}^{n_p}\textsf{T}_{j_1j_2\dots j_p\dots j_k}\textbf{v}_{j_p}, 
\end{equation}
which can be extended to 
\begin{equation}\label{eq6}
\begin{split}
\textsf{T}\times_1 \textbf{v}_1 \times_2\textbf{v}_2\times_3\dots \times_{k}\textbf{v}_k=\textsf{T}\textbf{v}_1\textbf{v}_2\dots\textbf{v}_k\in\mathbb{R}
\end{split}
\end{equation}
for $\textbf{v}_p\in \mathbb{R}^{ n_p}$. If $\textsf{T}$ is supersymmetric and $\textbf{v}_p=\textbf{v}$ for all $p$, the product (\ref{eq6}) is also known as the homogeneous polynomial associated with $\textsf{T}$ \cite{QI20051302,doi:10.1137/140989340}, and we write  it as $\textsf{T}\textbf{v}^{k}$ for simplicity. In other words, an $n$-dimensional homogeneous polynomial of degree $k$ can be uniquely determined by a $k$th-order $n$-dimensional supersymmetric tensor, analogous to quadratic forms (i.e., $\textbf{v}^\top\textbf{M}\textbf{v}$) in matrix theory. As an illustrative example, the two-dimensional homogeneous polynomial of degree three, i.e., 
$
f(x, y) = ax^3+bx^2y+cxy^2+dy^3
$,
can be represented by a supersymmetric tensor $\textsf{T}\in\mathbb{R}^{2\times 2\times 2}$ with entries  $\textsf{T}_{111}=a$, $\textsf{T}_{112}=\textsf{T}_{121}=\textsf{T}_{211}=\frac{b}{3}$, $\textsf{T}_{122}=\textsf{T}_{212}=\textsf{T}_{221}=\frac{c}{3}$, and $\textsf{T}_{222}=d$. Therefore, the product $\textsf{T}\textbf{v}^{k-1}\in\mathbb{R}^n$ belongs to the family of homogeneous polynomial systems (though it does not include the entire space for supersymmetric \textsf{T}, see Proposition 5).

\subsection{Tensor Eigenvalues}
The tensor eigenvalues of real supersymmetric tensors were first explored by Qi \cite{QI20051302, QI20071363} and Lim \cite{singularvaluetensor} independently. There are many notions of tensor eigenvalues including H-eigenvalues, Z-eigenvalues, M-eigenvalues, and U-eigenvalues \cite{chen2021multilinear,QI20051302, QI20071363}. Of particular interest of this paper are Z-eigenvalues (which are associated with tensor orthogonal decomposition). Given a $k$th-order supersymmetric tensor $\textsf{T}\in\mathbb{R}^{n\times n\times \dots \times n}$, the E-eigenvalues $\lambda\in\mathbb{C}$ and the E-eigenvectors $\textbf{v}\in\mathbb{C}^n$ of $\textsf{T}$ are defined as
\begin{equation}
\begin{cases}
   \textsf{T}\textbf{v}^{k-1} = \lambda\textbf{v}\\
   \textbf{v}^\top\textbf{v} = 1
\end{cases}.
\end{equation}
The E-eigenvalues $\lambda$ could be complex. If $\lambda$ are real, we call them Z-eigenvalues. It has been proved that a supersymmetric tensor always has Z-eigenvalues \cite{QI20071363}. Similar to matrix eigenvalues, the largest and smallest Z-eigenvalues of a supersymmetric tensor can be solved by
\begin{equation*}
    \max_{\textbf{v}\in\mathbb{R}^n}{\{\textsf{T}\textbf{v}^k:\|\textbf{v}\|_2=1\}} \text{ and } \min_{\textbf{v}\in\mathbb{R}^n}{\{\textsf{T}\textbf{v}^k:\|\textbf{v}\|_2=1\}},
\end{equation*}
respectively. Since the objective function is continuous and the feasible set is compact, the global maximizer and minimizer always exist. When $k=2$, the above optimization yields the matrix eigenvalue problem. Computing the E-eigenvalues or the Z-eigenvalues of a tensor is NP-hard \cite{Hillar:2013:MTP:2555516.2512329}. However, many numerical algorithms such as homotopy continuation approaches \cite{teneig,chen2016computing} and adaptive shifted power methods \cite{kolda2014adaptive} are proposed in order to best approximate the E-eigenvalues or the Z-eigenvalues of a tensor.

\begin{figure}[t]
    \centering
    \includegraphics[scale=0.32]{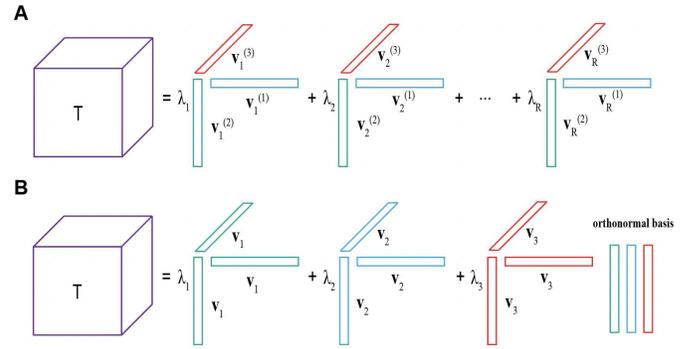}
    \caption{CP decomposition and orthogonal decomposition of a third-order three-dimensional tensor $\textsf{T}$. This figure is adapted from \cite{doi:10.1137/07070111X}.}
    \label{fig:1}
\end{figure}
\subsection{Tensor Orthogonal Decomposition}
Before talking about tensor orthogonal decomposition, we first introduce CANDECOMP/PARAFAC (CP) decomposition. Like rank-one matrices, a tensor $\textsf{T}\in\mathbb{R}^{n_1\times n_2\times \dots \times n_k}$ is rank-one if it can be written as the outer product of $k$ vectors, i.e.,
$
\textsf{T}=\textbf{v}^{(1)}\circ  \textbf{v}^{(2)} \circ \dots \circ \textbf{v}^{(k)}
$
(where ``$\circ$'' denotes the vector outer product). CP decomposition decomposes a tensor $\textsf{T}\in\mathbb{R}^{n_1\times n_2\times \dots \times n_k}$ into a sum of rank-one tensors, see Fig. 1 A. It is often useful to normalize all the vectors and have weights $\lambda_r$ in descending order in front:
\begin{equation}\label{eq1.7}
\textsf{T}=\sum_{r=1}^{R}\lambda_r\textbf{v}_r^{(1)}\circ  \textbf{v}_r^{(2)} \circ \dots \circ \textbf{v}_r^{(k)},
\end{equation}
where $\textbf{v}_r^{(p)}\in\mathbb{R}^{n_p}$ have unit length, and $R$ is called the CP rank of $\textsf{T}$ if it is the minimum integer that achieves (\ref{eq1.7}). Every tensor has a CP decomposition, and it is unique up to scaling and permutation under a weak condition, see details in \cite{doi:10.1137/07070111X}. The best CP rank approximation is ill-posed, but carefully truncating the CP rank will yield a good estimate of the original tensor \cite{chen2021multilinear}.

Tensor orthogonal decomposition is a special case of  CP decomposition. A $k$th-order supersymmetric tensor $\textsf{T}\in\mathbb{R}^{n\times n\times \dots \times n}$ is called orthogonally decomposable (odeco) if it can be written as a sum of vector outer products
\begin{equation}\label{eq:22}
\textsf{T} = \sum_{r=1}^n \lambda_r \textbf{v}_r\circ\textbf{v}_r\circ \stackrel{k}{\cdots} \circ\textbf{v}_r,
\end{equation}
where $\lambda_r\in\mathbb{R}$ in the descending order, and $\textbf{v}_r\in\mathbb{R}^n$ are orthonormal \cite{doi:10.1137/140989340}, see Fig. 1 B. It is easy to prove that $\lambda_r$ are the Z-eigenvalues of $\textsf{T}$ with the corresponding Z-eigenvectors $\textbf{v}_r$. Although $\lambda_r$ do not include all the Z-eigenvalues of $\textsf{T}$, it has been shown that the Z-spectral radius of \textsf{T} (i.e., the maximum absolute Z-eigenvalues) is equal to $\max{\{|\lambda_1|,|\lambda_n|\}}$ \cite{chen2022further}. In addition, not all supersymmetric tensors have orthogonal decomposition. Reobeva \cite{doi:10.1137/140989340} speculated that odeco tensors satisfy a set of polynomial equations that vanish on the odeco variety, which is the Zariski closure of the set of odeco tensors inside the space of $k$th-order $n$-dimensional complex supersymmetric tensors. Although the author only proved for the case when $n=2$, she provided strong evidence for its overall correctness. Furthermore, a tensor power method was proposed in \cite{doi:10.1137/140989340} in order to obtain the orthogonal decomposition of an odeco tensor.

\section{Main Results}\label{sec:3}
In this paper, we are interested in finding the explicit solution to a continuous-time homogeneous polynomial dynamical system (HPDS) of degree $k-1$ that can be represented by 
\begin{equation}\label{eq:2}
    \dot{\textbf{x}}(t)=\textsf{A}\times_1 \textbf{x}(t)\times_2\textbf{x}(t)\times_3\dots\times_{k-1}\textbf{x}(t)=\textsf{A}\textbf{x}(t)^{k-1},
\end{equation}
where $\textsf{A}\in\mathbb{R}^{n\times n\times \dots \times n}$ is a $k$th-order $n$-dimensional odeco tensor, and $\textbf{x}(t)\in\mathbb{R}^n$ is the state variable. We refer to such HPDS as odeco HPDS, see Section \ref{sec:6} for examples. 

\subsection{Explicit Solutions}
We find that the explicit solution to the odeco HPDS (\ref{eq:2}) can be solved in a simple fashion by exploiting the Z-eigenvalues and the Z-eigenvectors from the orthogonal decomposition of \textsf{A}. 

\textit{Proposition 1:} Suppose that $k\geq 3$ and $\textsf{A}\in\mathbb{R}^{n\times n\times \dots \times n}$ is odeco. Let the initial condition $\textbf{x}_0=\sum_{r=1}^n \alpha_r\textbf{v}_r$. Then the explicit solution to the odeco HPDS (\ref{eq:2}), given the initial condition $\textbf{x}_0$, is given by
\begin{equation}\label{eq:3}
    \textbf{x}(t) = \sum_{r=1}^n \Big (1-(k-2)\lambda_r\alpha_r^{k-2}t\Big)^{-\frac{1}{k-2}} \alpha_r\textbf{v}_r,
\end{equation}
where $\lambda_r$ are the Z-eigenvalues with the corresponding Z-eigenvectors $\textbf{v}_r$ from the orthogonal decomposition of $\textsf{A}$. Moreover, if $\lambda_r\alpha_r^{k-2}>0$ for some $r$, the solution (\ref{eq:3}) is only defined over the interval
\begin{equation}
    t\in \Big[0,  \min_{S}{\frac{1}{(k-2)\lambda_r\alpha_r^{k-2}}}\Big),
\end{equation}
where $S=\{r=1,2,\dots,n|\lambda_r\alpha_r^{k-2}>0\}$.

\textit{Proof:} Since $\textbf{v}_r$ are orthonormal, we can suppose that
\begin{equation*}
    \textbf{x}(t)=\sum_{r=1}^n c_r(t)\textbf{v}_r=\textbf{V}\textbf{c}(t),
\end{equation*}
where 
\begin{align*}
    \textbf{V} & =\begin{bmatrix} \textbf{v}_1 & \textbf{v}_2 & \dots & \textbf{v}_n
\end{bmatrix},\\
 \textbf{c}(t) & =\begin{bmatrix} c_1(t) & c_2(t) & \dots & c_n(t)\end{bmatrix}^\top.
\end{align*}
Clearly, $c_r(0) = \alpha_r$ for all $r$. Based on the property of tensor vector multiplications, it can be shown that
\begin{align*}
    \dot{\textbf{x}}(t) &= \Big(\sum_{r=1}^n \lambda_r \textbf{v}_r\circ\textbf{v}_r\circ \dots \circ\textbf{v}_r\Big)\times_1 \textbf{x}(t)\times_2\dots \times_{k-1}\textbf{x}(t)\\
    & = \Big(\sum_{r=1}^n \lambda_r \textbf{v}_r\circ\textbf{v}_r\circ \dots \circ\textbf{v}_r \Big)\times_1 \Big(\sum_{i=1}^n c_i(t)\textbf{v}_i \Big)\times_2 \dots \\ & \times_{k-1}\Big(\sum_{i=1}^n c_i(t)\textbf{v}_i\Big)
    = \sum_{r=1}^n \lambda_r \Big \langle \textbf{v}_r, \sum_{i=1}^n c_i(t)\textbf{v}_i \Big\rangle^{k-1} \textbf{v}_r\\&=\sum_{r=1}^n \lambda_rc_r(t)^{k-1}\textbf{v}_r.
\end{align*}
Thus, we have 
\begin{equation*}
    \dot{\textbf{c}}(t) = \boldsymbol{\lambda}* \textbf{c}(t)^{[k-1]} \Rightarrow  \dot{c}_r(t) = \lambda_r c_r(t)^{k-1},
\end{equation*}
where $\boldsymbol{\lambda}=\begin{bmatrix} \lambda_1 & \lambda_2 & \dots & \lambda_n\end{bmatrix}^\top$, $``*"$ denotes the element-wise multiplication, and the superscript $``[k-1]"$ denotes the element-wise $(k-1)$th power. By the method of separation of variables, we can solve for $c_r(t)$, which are given by
\begin{align*}
    &\int c_r(t)^{-(k-1)} dc_r(t) = \int \lambda_r dt \\ \Rightarrow &c_r(t) = \Big ((k-2)(w_r-\lambda_r t)\Big)^{-\frac{1}{k-2}}.
\end{align*}
Thus, plugging the initial condition yields
\begin{equation*}
    c_r(t) = \Big (1-(k-2)\lambda_r\alpha_r^{k-2}t\Big)^{-\frac{1}{k-2}} \alpha_r,
\end{equation*}
and the result follows immediately. Moreover, if  $\lambda_r\alpha_r^{k-2}>0$ for some $r$, the corresponding coefficient functions $c_r(t)$ will have singularities at $t=\frac{1}{(k-2)\lambda_r\alpha_r^{k-2}}$. Thus, the domains of $c_r(t)$ are given by $t\in\Big [0, \frac{1}{(k-2)\lambda_r\alpha_r^{k-2}}\Big)$. The other branches of $c_r(t)$ over $t\in \Big(\frac{1}{(k-2)\lambda_r\alpha_r^{k-2}}, \infty\Big)$ do not satisfy the initial conditions, so they are not included in the solutions of $c_r(t)$. Therefore, the overall domain of the solution (\ref{eq:3}) is given by
\begin{equation*}
    D = \bigcap_S \Big [0, \frac{1}{(k-2)\lambda_r\alpha_r^{k-2}}\Big)=\Big[0,  \min_{S}{\frac{1}{(k-2)\lambda_r\alpha_r^{k-2}}}\Big),
\end{equation*}
where $S=\{r=1,2,\dots,n|\lambda_r\alpha_r^{k-2}>0\}$. Note that if $\lambda_r\alpha_r^{k-2}\leq 0$ for all $r$, the domain of the solution (\ref{eq:3}) will be $D=[0, \infty)$. \hfill $\blacksquare$

The coefficients $\alpha_r$ can be found from the inner product between $\textbf{x}_0$ and $\textbf{v}_r$. Moreover, when $\lambda_r>0$ for some $r$, the solution (\ref{eq:3}) has singularity points, and the system blows up within a finite time. In general, analyzing the occurrence of blow-up of solutions to polynomial systems is challenging \cite{goriely1998finite}. Thus, our approach offers a simple way to detect the singularity behaviors of a HPDS. When $k=2$, the result reduces to the linear systems' solutions, i.e., 
\begin{align*}
    \lim_{k \rightarrow 2}\textbf{x}(t)&=\lim_{k \rightarrow 2} \sum_{r=1}^n \Big (1-(k-2)\lambda_r\alpha_r^{k-2}t\Big)^{-\frac{1}{k-2}} \alpha_r\textbf{v}_r \\&= \lim_{p\rightarrow \infty} \sum_{r=1}^n\Big(1+\frac{\lambda_rt}{p}\Big)^p\alpha_r\textbf{v}_r = \sum_{r=1}^n \exp{\{\lambda_rt\}}\alpha_r\textbf{v}_r,
\end{align*}
where $\lambda_r$ become the eigenvalues of the dynamic matrix with the corresponding eigenvectors $\textbf{v}_r$.

\subsection{Stability}
In linear control theory, it is conventional to investigate so-called (internal) stability \cite{rugh1996linear}. We are able to discuss the stability properties of an odeco HPDS based on the explicit solution (\ref{eq:3}). First, we explore the number of equilibrium points of an odeco HPDS, which is similar to the cases in linear systems.

\textit{Proposition 2:} The odeco HPDS (\ref{eq:2}) has a unique equilibrium point at the origin if $\lambda_r\neq 0$ for all $r=1,2,\dots,n$, where $\lambda_r$ are the Z-eigenvalues from the orthogonal decomposition of \textsf{A}. Otherwise, it has infinitely many equilibrium points.

\textit{Proof:} Suppose that the equilibrium point $\textbf{x}_e=\sum_{r=1}^n e_r\textbf{v}_r$ (where $\textbf{v}_r$ are the Z-eigenvectors). Similarly, it can be shown that
\begin{equation*}
    \textsf{A}\textbf{x}_e^{k-1} = \sum_{r=1}^n\lambda_r e_r^{k-1}\textbf{v}_r = \textbf{0}.
\end{equation*}
Since $\textbf{v}_r$ are linearly independent, $\lambda_r e_r^{k-1}=0$. Thus, if $\lambda_r\neq 0$, then $e_r=0$ for all $r$, which implies that the odeco HPDS (\ref{eq:2}) has a unique equilibrium point at the origin. If $\lambda_r = 0$ for some $r$, the corresponding $e_r$ can be chosen arbitrarily, and the system has infinitely many equilibrium points. \hfill $\blacksquare$

We only need to focus our attention on the equilibrium point at the origin since the behaviors of other equilibrium points will be the same (similar to linear systems). The equilibrium point $\textbf{x}_e=\textbf{0}$ of an odeco HPDS is called stable if $\|\textbf{x}(t)\|\leq \gamma \|\textbf{x}_0\|$ for some initial condition $\textbf{x}_0$ and $\gamma>0$, asymptotically stable if $\|\textbf{x}(t)\|\rightarrow 0$ as $t\rightarrow \infty$, and unstable if  $\|\textbf{x}(t)\|\rightarrow \infty$ as $t\rightarrow c$ for $c>0$ (which is also referred to as finite-time blow-up in the literature \cite{goriely1998finite}). Here $``\|\cdot\|"$ denotes the Frobenius norm. We discover that the stability properties of the odeco HPDS (\ref{eq:2}) at the equilibrium point $\textbf{x}_e=\textbf{0}$ are similar to those of linear systems, but depend on both the Z-eigenvalues of $\textsf{A}$ and initial conditions. 

\textit{Proposition 3:} 
Suppose that $k\geq 3$. Let the initial condition $\textbf{x}_0=\sum_{r=1}^n \alpha_r\textbf{v}_r$. For the odeco HPDS (\ref{eq:2}), the equilibrium point $\textbf{x}_e=\textbf{0}$ is:
\begin{enumerate}
    \item stable if and only if $\lambda_r\alpha_r^{k-2}\leq 0$ for all $r=1,2,\dots,n$;
    \item asymptotically stable if and only if $\lambda_r\alpha_r^{k-2}<0$ for all $r=1,2,\dots,n$;
    \item unstable if and only if $\lambda_r\alpha_r^{k-2}>0$ for some $r=1,2,\dots,n$, 
\end{enumerate}
where $\lambda_r$ are the Z-eigenvalues from the orthogonal decomposition of $\textsf{A}$.

\textit{Proof:} By the triangle inequality, it can be shown that 
\begin{equation*}
    \|\textbf{x}(t)\| = \|\sum_{r=1}^nc_r(t)\textbf{v}_r\|\leq \sum_{r=1}^n|c_r(t)|\|\textbf{v}_r\| = \sum_{r=1}^n|c_r(t)|.
\end{equation*}
Since $\lambda_r\alpha_r^{k-2}\leq 0$ for all $r=1,2,\dots,n$, the coefficient functions $|c_r(t)|$ are bounded by $|\alpha_r|$ over $t\in[0,\infty)$. Then we have
\begin{equation*}
    \|\textbf{x}(t)\|\leq \sum_{r=1}^n|\alpha_r| =\|\textbf{x}_0\|_1\leq \sqrt{n}\|\textbf{x}_0\|.
\end{equation*}
Therefore, the equilibrium point $\textbf{x}_e=\textbf{0}$ is stable. On the other hand, since $\textbf{v}_r$ are orthonormal,  $\|\textbf{x}(t)\|=\|\textbf{V}\textbf{c}(t)\|=\|\textbf{c}(t)\|$ where $\textbf{V}$ and $\textbf{c}(t)$ are the same as defined in Proposition 1. If $\|\textbf{x}(t)\|=\|\textbf{c}(t)\|\leq \gamma \|\textbf{x}_0\|$, all the coefficient functions $c_r(t)$ must be bounded for $t\geq 0$. Thus, $\lambda_r\alpha_r^{k-2}$ must lie in the closed left-half plane for all $r$. The other two cases can be shown similarly. \hfill $\blacksquare$

When $k=2$, the above conditions reduce to the famous linear stability conditions. The inequalities obtained from the asymptotic stability condition  can provide us with the region of attraction of the odeco HPDS (\ref{eq:2}), i.e., 
\begin{equation}
    R = \{\textbf{x}: \lambda_r\alpha_r^{k-2} <0  \text{ where }\textbf{x}=\sum_{r=1}^n\alpha_r\textbf{v}_r\},
\end{equation}
where $\textbf{v}_r$ are the Z-eigenvectors from the orthogonal decomposition of $\textsf{A}$ corresponding to the Z-eigenvalues $\lambda_r$. Furthermore, when $k$ is even, $\alpha_r^{k-2}$ will always be  greater than or equal to zero. Thus, the stability conditions can be simplified for the odeco HPDS (\ref{eq:2}) of odd degree. 

\textit{Corollary 1:} Suppose that $k\geq 4$ is even. For the odeco HPDS (\ref{eq:2}), the equilibrium point $\textbf{x}_e=\textbf{0}$ is:
\begin{enumerate}
    \item stable if and only if $\lambda_r\leq 0$ for all $r=1,2,\dots,n$;
    \item asymptotically stable if and only if $\lambda_r<0$ for all $r=1,2,\dots,n$;
    \item unstable if and only if $\lambda_r>0$ for some $r=1,2,\dots,n$,
\end{enumerate}
where $\lambda_r$ are the Z-eigenvalues from the orthogonal decomposition of $\textsf{A}$.

\textit{Proof:} The results follow immediately from Proposition 3 when $k$ is even. \hfill $\blacksquare$

When $k$ is even, the stability conditions are exactly the same as those from linear systems, i.e., the odeco HPDS (\ref{eq:2}) of odd degree is globally stable if and only if all the Z-eigenvalues $\lambda_r$ from the orthogonal decomposition of $\textsf{A}$ lie in the left-half plane. However, as mentioned, computing the Z-eigenvalues of a supersymmetric tensor is NP-hard  \cite{Hillar:2013:MTP:2555516.2512329}. If we know an  upper bound of the largest Z-eigenvalue of a supersymmetric tensor, it will save a great amount of computation to determine the asymptotic stability of the odeco HPDS (\ref{eq:2}). Chen \cite{chen2021stability} found that the largest Z-eigenvalue of an even-order supersymmetric tensor is upper bounded by the largest eigenvalue of one of its unfolded matrices. 

\textit{Lemma 1:} 
Let $\textsf{A}\in\mathbb{R}^{n\times n\times \dots \times n}$ be an even-order supersymmetric tensor. Then the largest Z-eigenvalue $\lambda_{\max}$ of $\textsf{A}$ is upper bounded by $\mu_{\max}$ where $\mu_{\max}$ is the largest eigenvalue of $\psi(\textsf{A})$ defined as
\begin{equation}\label{eq:4}
    \textbf{A} = \psi(\textsf{A}) \text{ such that } \textsf{A}_{j_1i_1\dots j_ki_k} \xrightarrow{\psi} \textbf{A}_{ji},
\end{equation}
with $j=j_1+\sum_{p=2}^k(j_p-1)n^{p-1}$ and $i=i_1+\sum_{p=2}^k(i_p-1)n^{p-1}$.

\textit{Corollary 2:}
Suppose that $k\geq 4$ is even. For the odeco HPDS (\ref{eq:2}), the equilibrium point $\textbf{x}=\textbf{0}$ is:
\begin{enumerate}
    \item stable if $\mu_{\max} \leq 0$;
    \item asymptotically stable if $\mu_{\max}<0$, 
\end{enumerate}
where $\mu_{\max}$ is the largest eigenvalue of $\psi(\textsf{A})$ defined in (\ref{eq:4}).

\textit{Proof:} Based on Lemma 1, we know that $\lambda_1\leq \lambda_{\max}\leq \mu_{\max}$. Therefore,  the result follows immediately from Corollary 1.  \hfill $\blacksquare$

Note that $\lambda_1$ is the largest Z-eigenvalue from the orthogonal decomposition of $\textsf{A}$, while $\lambda_{\max}$ is the largest Z-eigenvalue of $\textsf{A}$. There are many other upper bounds for the largest Z-eigenvalue of a supersymmetric tensor  \cite{chang2013some,he2014upper,ma2019some,wu2018upper}. Given an odeco dynamic tensor, the better upper bound of the largest Z-eigenvalue, the more strong stability conditions we can obtain.

\section{Odeco HPDS with Constant Control}\label{sec:4}
The constant control problem has arisen in many dynamical systems and control applications. For example,  model predictive control has been successfully implemented for stable plants based on linear models by optimizing a constant input on the whole horizon \cite{qin1997overview}. In the context of population dynamics modeled by HPDS, constant inputs can be viewed as migration or supply rates for different population groups. Therefore, it is significant to investigate the system properties of an odeco HPDS with constant control, i.e., 
\begin{equation}\label{eq:5}
    \dot{\textbf{x}}(t) = \textsf{A}\textbf{x}(t)^{k-1} + \textbf{b},
\end{equation}
where $\textsf{A}\in\mathbb{R}^{n\times n\times \dots \times n}$ is a $k$th-order $n$-dimensional odeco tensor, and $\textbf{b}\in\mathbb{R}^n$ is a constant control vector. We find that the complete solution to the polynomial dynamical system (\ref{eq:5}) can be solved implicitly by using Gauss hypergeometric functions.

\textit{Proposition 4:} Suppose that $k\geq 3$. Let $\textbf{x}(t)=\sum_{r=1}^n c_r(t)\textbf{v}_r$ with the initial conditions $c_r(0)=\alpha_r$. For the odeco HPDS with constant control (\ref{eq:5}), the coefficient functions $c_r(t)$ can be solved implicitly by 
\begin{equation}\label{eq:7}
    t=-\frac{g\Big(\frac{k-2}{k-1}, -\frac{\tilde{b}_r}{\lambda_rc_r(t)^{k-1}}\Big)}{(k-2)\lambda_rc_r(t)^{k-2}} +\frac{g\Big(\frac{k-2}{k-1}, -\frac{\tilde{b}_r}{\lambda_r\alpha_r^{k-1}}\Big)}{(k-2)\lambda_r\alpha_r^{k-2}},
\end{equation}
where $\lambda_r$ are the Z-eigenvalues with the corresponding Z-eigenvectors $\textbf{v}_r$ from the orthogonal decomposition of $\textsf{A}$, $\tilde{b}_r$ are the $r$th entries of $\textbf{V}^\top\textbf{b}$ (\textbf{V} contains all the vectors $\textbf{v}_r$), and $g(\cdot, \cdot)$ is the specified Gauss hypergeometric function \cite{jan2012use} defined as
\begin{equation*}
    g(a, z) = {}_2 F_1(1, a; a+1;z)=a\sum_{m=0}^{\infty} \frac{z^m}{a+m}.
\end{equation*}

\textit{Proof:}
Since $\textbf{x}(t)=\sum_{r=1}^n c_r(t)\textbf{v}_r$, we can rewrite the polynomial dynamical system (\ref{eq:5}) as follows:
\begin{align*}
    \textbf{V}\dot{\textbf{c}}(t) &=\textbf{V} (\boldsymbol{\lambda}*\textbf{c}(t)^{[k-1]})+\textbf{V}\textbf{V}^\top \textbf{b}\\
    \Rightarrow \dot{\textbf{c}}(t) &= \boldsymbol{\lambda}*\textbf{c}(t)^{[k-1]}+\tilde{\textbf{b}},
\end{align*}
where $\tilde{\textbf{b}}=\textbf{V}^\top\textbf{b}$. Therefore, for each coefficient function $c_r(t)$, we have
\begin{equation}\label{eq:6}
    \dot{c}_r(t) = \lambda_r c_r(t)^{k-1}+\tilde{b}_r.
\end{equation}
The differential equation (\ref{eq:6}) is a particular form of the Chini's equation \cite{kamke2013differentialgleichungen}, and can be solved implicitly by using Gauss hypergeometric functions. Based on the method of separation of variables, it can be shown that 
\begin{align*}
    &\int \frac{1}{\lambda_rc_r(t)^{k-1}+\tilde{b}_r} dc_r(t) = \int 1 dt \\ \Rightarrow &-\frac{g\Big(\frac{k-2}{k-1}, -\frac{\tilde{b}_r}{\lambda_rc_r(t)^{k-1}}\Big)}{(k-2)\lambda_rc_r(t)^{k-2}} = t+w_r.
\end{align*}
Plugging the initial conditions yields
\begin{equation*}
    t=-\frac{g\Big(\frac{k-2}{k-1}, -\frac{\tilde{b}_r}{\lambda_rc_r(t)^{k-1}}\Big)}{(k-2)\lambda_rc_r(t)^{k-2}} +\frac{g\Big(\frac{k-2}{k-1}, -\frac{\tilde{b}_r}{\lambda_r\alpha_r^{k-1}}\Big)}{(k-2)\lambda_r\alpha_r^{k-2}},
\end{equation*}
and the proof is complete. \hfill $\blacksquare$

The solutions of $c_r(t)$ can be further solved by any nonlinear  solver given a specific time point $t$. We then can recover the complete solution of $\textbf{x}(t)$ based on the values of $c_r(t)$. Moreover, although $g(a, z)$ is defined for $|z|<1$, it can be  analytically continued along any path in the complex plane that avoids the branch points one and infinity \cite{gasper2004basic}. When $k=3$, the differential equation (\ref{eq:6}) is also known as the Riccati equation, which can be converted to a second-order linear system.


\section{Extension to General HPDS}\label{sec:5}
In this section, we  extend the previous results to general HPDS. First, we introduce the notion of almost symmetric for cubical tensors. 

\textit{Definition 1:} A $k$th-order $n$-dimensional tensor $\textsf{A}\in\mathbb{R}^{n\times n\times \dots \times n}$ is called almost symmetric if it is symmetric only respect to its first $k-1$ modes. 

\textit{Proposition 5:}
Every HPDS of degree $k-1$ can be represented by
\begin{equation}\label{eq:10}
    \dot{\textbf{x}}(t) = \textsf{A}\textbf{x}(t)^{k-1},
\end{equation}
where $\textsf{A}\in\mathbb{R}^{n\times n \times \dots \times n}$ is a $k$th-order almost symmetric tensor.

\textit{Proof:} Since $\textsf{A}$ is almost symmetric, its $(k-1)$th-order sub-tensors $\textsf{A}_{::\dots j}$ are supersymmetric for $j=1,2,\dots,n$. We know that every homogeneous polynomial of degree $k-1$ can be uniquely represented by a $(k-1)$th-order supersymmetric tensor. Therefore, the result follows immediately.  \hfill $\blacksquare$

The colon operator ``:" in the proof acts as shorthand to include all subscripts in a particular array dimension as used in MATLAB. Since $\textsf{A}$ is almost symmetric, a CP decomposition of \textsf{A} is given by
\begin{equation}\label{eq:9}
    \textsf{A} = \sum_{r=1}^R  \textbf{v}_r\circ \textbf{v}_r \circ \dots \circ\textbf{v}^{(f)}_r.
\end{equation}
Without loss of generality, we multiply the weights $\lambda_r$ into the vector $\textbf{v}^{(f)}_r$ beforehand. Our goal is to construct a linear transformation $\textbf{P}\in\mathbb{R}^{n\times n}$ with $\textbf{x}(t) = \textbf{P}\textbf{y}(t)$ such that the transformed system can be represented by an odeco tensor, i.e., 
\begin{equation}\label{eq:11}
    \dot{\textbf{y}}(t) = \tilde{\textsf{A}}\textbf{y}(t)^{k-1},
\end{equation}
where $\tilde{\textsf{A}}\in\mathbb{R}^{n\times n \times \dots \times n}$ is a $k$th-order odeco tensor. 

\textit{Proposition 6:}
Suppose that $\textsf{A}\in\mathbb{R}^{n\times n \times \dots \times n}$ has the CP decomposition (\ref{eq:9}) with $R=n$. If there exist an invertible linear transformation $\textbf{P}\in\mathbb{R}^{n\times n}$ and a diagonal matrix $\boldsymbol{\Lambda}\in\mathbb{R}^{n\times n}$ that satisfy the following conditions:
\begin{enumerate}
    \item $\textbf{P}^\top\textbf{V} = \textbf{P}^{-1}\textbf{V}^{(f)}\boldsymbol{\Lambda}^{-1}$;
    \item $\textbf{P}^\top\textbf{V}$ is an orthogonal matrix,
\end{enumerate}
where $\textbf{V}\in\mathbb{R}^{n\times n}$ and $\textbf{V}^{(f)}\in\mathbb{R}^{n\times n}$ are the matrices that contain all the vectors $\textbf{v}_r$ and $\textbf{v}_r^{(f)}$, respectively, then the HPDS (\ref{eq:10}) can be transformed to the odeco HPDS (\ref{eq:11}).

\begin{algorithm}[t]
\caption{Determining if a general HPDS can be transformed to an odeco HPDS.}
\label{alg:1}
\begin{algorithmic}[1]
\STATE{Given a general HPDS of the form (\ref{eq:10}) and a threshold $\epsilon$ (default: $\epsilon=10^{-14}$)}\\
\STATE{Create a CP decomposition model\\\centerline{\texttt{model = struct}}}
\STATE{Randomly initialize a variable $\textbf{R}\in\mathbb{R}^{n\times n}$ and a weight vector $\boldsymbol{\lambda}\in\mathbb{R}^n$\\ \centerline{\texttt{model.variable.R = randn(n,n)}}\centerline{\texttt{model.variable.w = randn(1,n)}}}
\STATE{Impose the conditions on the structure of the factor matrices such that $\textbf{V}=\textbf{R}$ and $\textbf{V}^{(f)}=(\textbf{R}^{-1})^\top$ \\ \centerline{\texttt{model.factor.w = `w'}}\centerline{\texttt{model.factor.V = `R'}} \centerline{\texttt{model.factor.Vf = \{`R',@(z,task)}}\centerline{\texttt{struct\_invtransp(z, task)\}}}}
\STATE{Compute the CP decomposition of the dynamic tensor \textsf{A} based on the imposed conditions\\ \centerline{\texttt{model.factorizations.symm.data = A}}\centerline{\texttt{model.factorizations.cpd =}} \centerline{\texttt{\{`V',$\dots$,`V',`Vf',`w'\}}} \centerline{\texttt{cpd = sdf\_nls(model)}}}
\STATE{Use the obtained factor matrices $\textbf{V}$ and $\textbf{V}^{(f)}$ with the weight vector $\boldsymbol{\lambda}$ to build the estimated tensor $\hat{\textsf{A}}$\\\centerline{\texttt{w = cpd.factors.w}}\centerline{\texttt{V = cpd.factors.V}} \centerline{\texttt{Vf = cpd.factors.Vf}}\centerline{\texttt{Ahat = cpdgen(\{V,$\dots$,V,Vf,w\})}}}
\IF{$\|\textsf{A}-\hat{\textsf{A}}\|< \epsilon$}
\STATE{The HPDS with dynamic tensor \textsf{A} can be transformed to an odeco HPDS}\\
\ENDIF
\end{algorithmic}
\end{algorithm}

\textit{Proof:}
Since $\textbf{y}(t) = \textbf{P}^{-1}\textbf{x}(t)$, we can write 
\begin{align*}
    \dot{\textbf{y}}(t) &= \textbf{P}^{-1}\dot{\textbf{x}}(t) = \textbf{P}^{-1}\Big (\textsf{A}\textbf{x}(t)^{k-1}\Big)\\
    &= \textbf{P}^{-1}\Big[\Big(\sum_{r=1}^n  \textbf{v}_r\circ \textbf{v}_r \circ \dots \circ\textbf{v}^{(f)}_r\Big) \Big(\textbf{P}\textbf{y}(t)\Big)^{k-1}\Big]\\
    & = \Big(\sum_{r=1}^n  \textbf{P}^\top\textbf{v}_r\circ \textbf{P}^\top\textbf{v}_r \circ \dots \circ\textbf{P}^{-1}\textbf{v}^{(f)}_r\Big)\textbf{y}(t)^{k-1}.
\end{align*}
If $\textbf{P}^\top\textbf{V} = \textbf{P}^{-1}\textbf{V}^{(f)}\boldsymbol{\Lambda}^{-1}$, then 
\begin{align*}
     &\Big(\sum_{r=1}^n  \textbf{P}^\top\textbf{v}_r\circ \textbf{P}^\top\textbf{v}_r \circ \dots \circ\textbf{P}^{-1}\textbf{v}^{(f)}_r\Big)\textbf{y}(t)^{k-1}\\ =& \Big(\sum_{r=1}^n  \textbf{P}^\top\textbf{v}_r\circ \textbf{P}^\top\textbf{v}_r \circ \dots \circ\lambda_r\textbf{P}^\top\textbf{v}_r\Big)\textbf{y}(t)^{k-1},
\end{align*}
where $\lambda_r$ are the $r$th diagonal of $\boldsymbol{\Lambda}$. Moreover, if $\textbf{P}^\top\textbf{V}$ is an orthogonal matrix,  the transformed dynamic tensor 
\begin{align*}
    \tilde{\textsf{A}} = \sum_{r=1}^n \lambda_r \textbf{P}^\top\textbf{v}_r\circ \textbf{P}^\top\textbf{v}_r \circ \dots \circ\textbf{P}^\top\textbf{v}_r
\end{align*}
is odeco. Thus, the result follows immediately. \hfill $\blacksquare$

The above two conditions can be further simplified.

\textit{Corollary 3:} Suppose that $\textsf{A}\in\mathbb{R}^{n\times n \times \dots \times n}$ has the CP decomposition (\ref{eq:9}) with $R=n$. Let $\textbf{V}\in\mathbb{R}^{n\times n}$ and $\textbf{V}^{(f)}\in\mathbb{R}^{n\times n}$ be the matrices that contain all the vectors $\textbf{v}_r$ and $\textbf{v}_r^{(f)}$, respectively. Let $\textbf{W}=(\textbf{V}^{-1})^\top$. If there exist $\lambda_r\in\mathbb{R}$ such that $\textbf{w}_r=\lambda_r^{-1}\textbf{v}^{(f)}_r$ for all $r$ ($\textbf{w}_r$ are the column vectors of $\textbf{W}$), then the HPDS (\ref{eq:10}) can be transformed to the odeco HPDS (\ref{eq:11}).

\textit{Proof:} The result follows immediately by combining the two conditions from Proposition 6, i.e., 
$
    \textbf{P}\textbf{P}^\top\textbf{V} = \textbf{V}^{(f)}\boldsymbol{\Lambda}^{-1}\Rightarrow \textbf{P}\textbf{P}^{-1}\textbf{W} = \textbf{V}^{(f)}\boldsymbol{\Lambda}^{-1} \Rightarrow \textbf{W} = \textbf{V}^{(f)}\boldsymbol{\Lambda}^{-1}.
$ 
\hfill $\blacksquare$

\begin{figure}[t]
    \centering
    \includegraphics[scale=0.2]{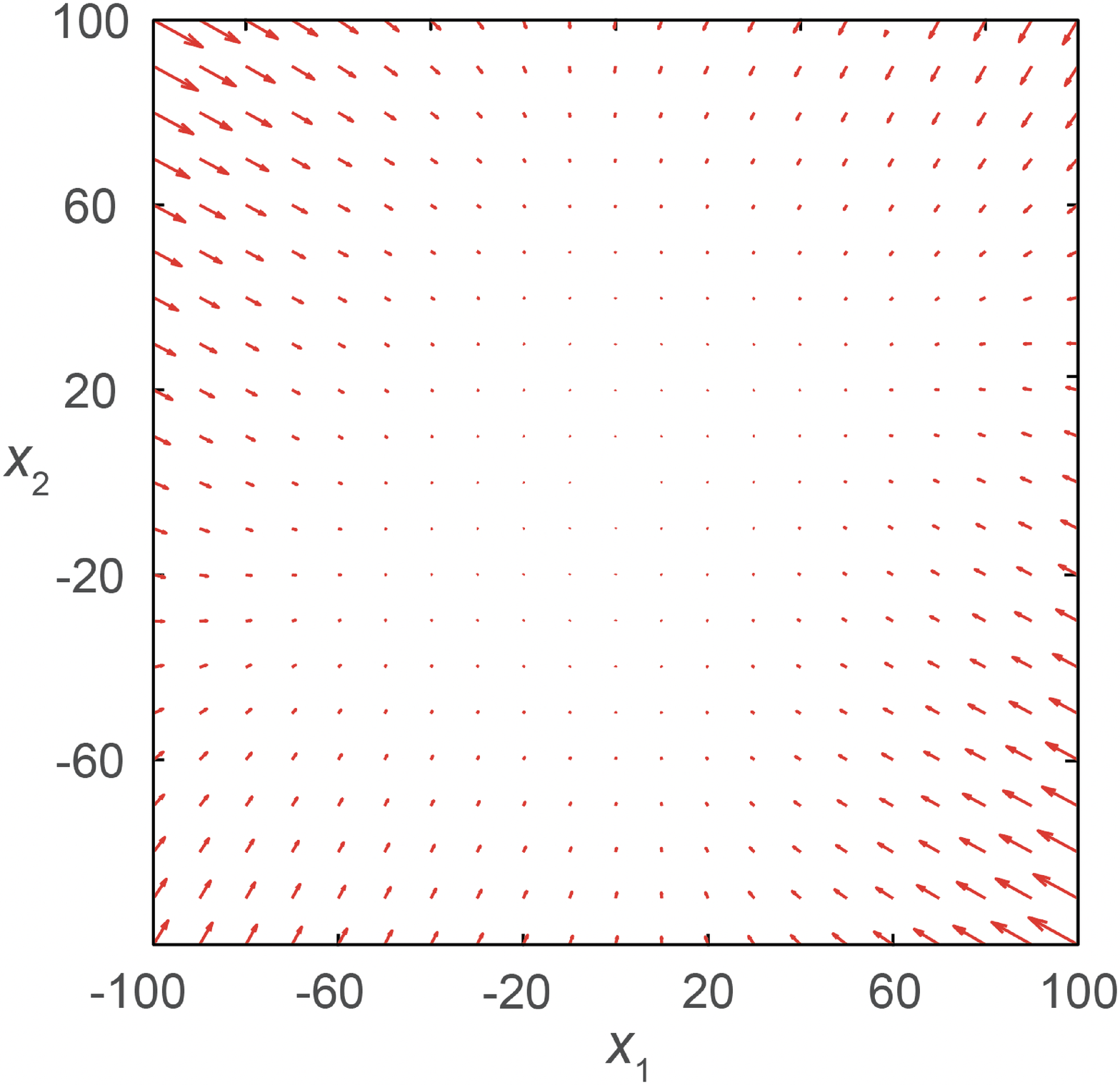}
    \caption{Synthetic example. Vector field plot of the odeco HPDS.}
    \label{fig:2}
\end{figure}

The MATLAB toolbox TensorLab \cite{vervliet2016tensorlab} offers a nice feature in computing a CP decomposition of a tensor by imposing specific structure on the factor matrices. Detailed steps of determining if a general HPDS can be converted to an odeco HPDS are summarized in Algorithm 1, where ``$\|\cdot\|$'' denotes the tensor Frobenius norm in Step 7, and all the syntaxes can be found in the TensorLab toolbox. Although the conditions imposed on the factor matrices are restricted, one can adjust the values of $\epsilon$ to obtain an approximated odeco HPDS regardless of the CP rank of the dynamic tensor. In other words, Algorithm 1 can be applied to any dynamic tensor and return an approximated odeco tensor with an error up to threshold $\epsilon$. As long as $\epsilon$ is selected reasonably, the properties of the approximated odeco HPDS should be close to these of the original HPDS. After determining that a HPDS can be converted to or approximated by an odeco HPDS, one can compute the linear transformation $\textbf{P}$ for an arbitrary orthonormal basis, i.e., $\textbf{P}^\top\textbf{V} = \textbf{U}$, where $\textbf{U}$ is an arbitrary orthogonal matrix.  Since the transformed state $\textbf{y}(t)$ can be solved explicitly, the solution of the original HPDS is then given by $\textbf{x}(t)=\textbf{P}\textbf{y}(t)$. The results of stability also follow.

\section{Numerical Examples} \label{sec:6}
All the numerical examples presented were performed on a Macintosh machine with 16 GB RAM and a 2 GHz Quad-Core Intel Core i5 processor in MATLAB R2020b.

\subsection{Synthetic Example}
Given a following two-dimensional odeco HPDS of degree three
\begin{equation*}
    \begin{cases}
    \dot{x}_1 &= - 1.2593x_1^3 + 1.6630x_1^2x_2 - 1.5554x_1x_2^2 - 0.1386x_2^3\\
    \dot{x}_2 & = 0.5543x_1^3 - 1.5554x_1^2x_2 - 0.4158x_1x_2^2 - 0.7037x_2^3\\
    \end{cases},
\end{equation*}
it can be represented in the form of (\ref{eq:2}) with
\begin{align*}
\textsf{A}_{::11} &= \begin{bmatrix} -1.2593 &   0.5543\\
    0.5543 &  -0.5185\end{bmatrix} \text{ } \textsf{A}_{::12} = \begin{bmatrix} 0.5543 &  -0.5185\\
   -0.5185 &  -0.1386\end{bmatrix}\\
\textsf{A}_{::21} &= \begin{bmatrix} 0.5543 &  -0.5185\\
   -0.5185 &  -0.1386\end{bmatrix} \text{ } \textsf{A}_{::22} = \begin{bmatrix} -0.5185 &  -0.1386\\
   -0.1386 &  -0.7037\end{bmatrix}
\end{align*}
such that $\textsf{A}$ is odeco. The two Z-eigenvalues in the orthogonal decomposition of $\textsf{A}$ are $\lambda_1=-1$ and $\lambda_2=-2$. Therefore, according to Corollary 1, the odeco HPDS is globally asymptotically stable. The vector field of the system is shown in Fig. 2. 

\begin{figure}[t]
    \centering
    \includegraphics[scale=0.29]{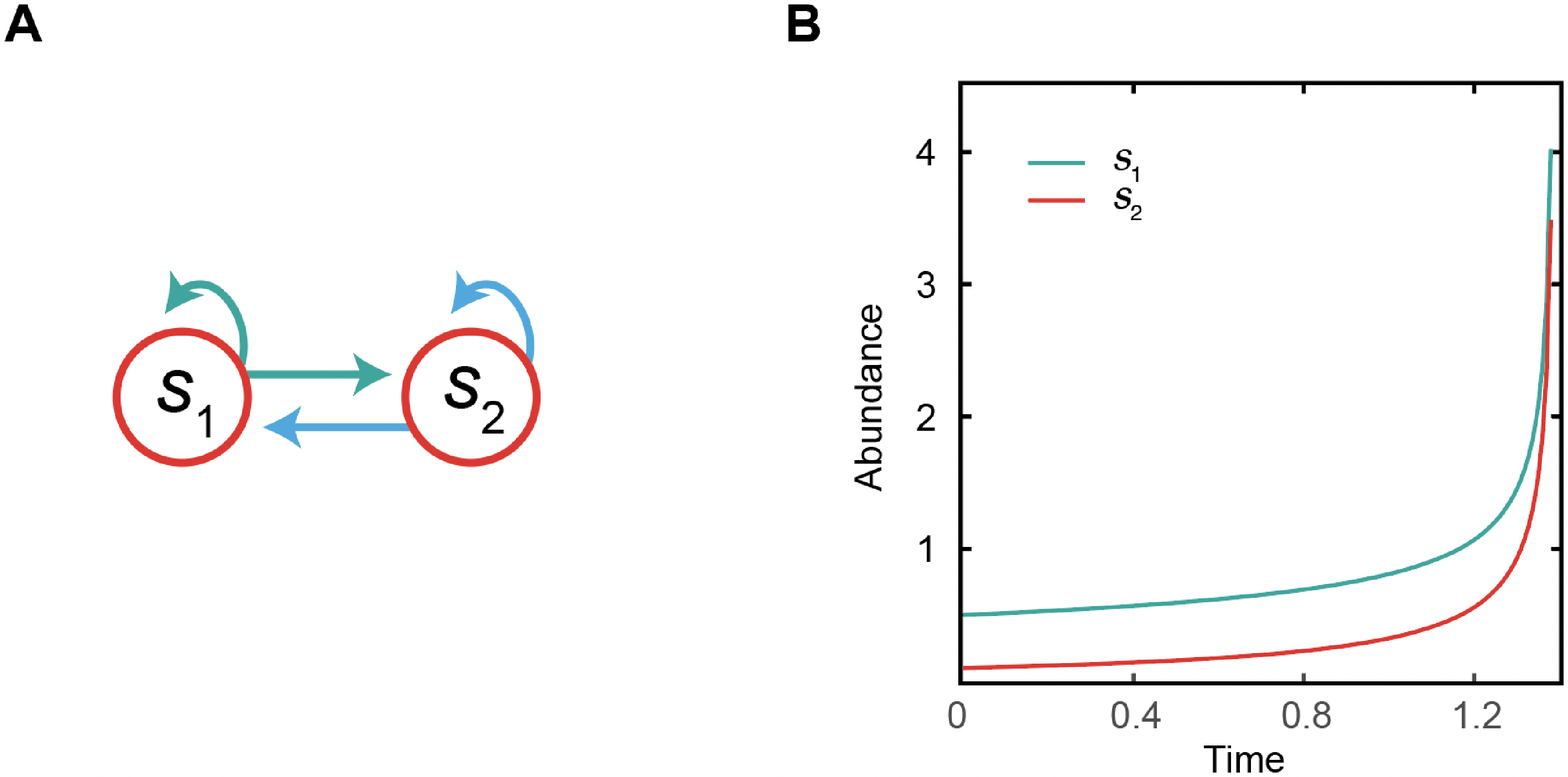}
    \caption{Population model example. (A) Schematic plot of two species with interactions described by the odeco HPDS. (B) Plot of the trajectories of the population dynamics. We used the initial condition $\textbf{x}_0=\begin{bmatrix}0.5 & 0.1\end{bmatrix}^\top$ in generating the plot.}
    \label{fig:3}
\end{figure}

\subsection{Population Model}
The goal of this example is to show the odeco structure in modeling population dynamics in ecological systems. Suppose that there are two species $s_1$ and $s_2$ with state variables $x_1$ and $x_2$ representing the abundance of the two species, the two factor vectors (which forms an orthonormal basis) for constructing the odeco dynamic tensor are given by
\begin{equation*}
    \textbf{v}_1=\begin{bmatrix}\frac{\sqrt{2}}{2}\\\frac{\sqrt{2}}{2}\end{bmatrix} \text{ and } \textbf{v}_2=\begin{bmatrix}\frac{\sqrt{2}}{2}\\-\frac{\sqrt{2}}{2}\end{bmatrix},
\end{equation*}
and the order of the tensor is four (i.e., $k=4$). Each factor vector is able to tell the latent interactions between the two species. The first vector $\textbf{v}_1$ indicates that the two species are self-promoted with mutual positive regulations, while the second vector $\textbf{v}_2$ implies that the two species are self-promoted with mutual negative regulations. The interactions between the two species become clearer when we have the expanded forms generated by $\textbf{v}_1$ and $\textbf{v}_2$, respectively, i.e.,
\begin{align*}
    \begin{cases}
    \dot{x}_1 & = \frac{1}{4}x_1^3 + (\frac{3}{4}x_1^2x_2+\frac{3}{4}x_1x_2^2+\frac{1}{4}x_2^3)\\
    \dot{x}_2 & = \frac{1}{4}x_2^3 + (\frac{3}{4}x_2^2x_1+\frac{3}{4}x_2x_1^2+\frac{1}{4}x_1^3)
    \end{cases},\\
    \begin{cases}
    \dot{x}_1 & = \frac{1}{4}x_1^3 - (\frac{3}{4}x_1^2x_2-\frac{3}{4}x_1x_2^2+\frac{1}{4}x_2^3)\\
    \dot{x}_2 & = \frac{1}{4}x_2^3 - (\frac{3}{4}x_2^2x_1-\frac{3}{4}x_2x_1^2+\frac{1}{4}x_1^3)
    \end{cases}.
\end{align*}
Suppose that $\lambda_1=\lambda_2=2$. Therefore, the overall population dynamics between the two species is given by
\begin{equation*}
    \begin{cases}
    \dot{x}_1 & = x_1^3 + 3x_1x_2^2\\
    \dot{x}_2 & = x_2^3 + 3x_2x_1^2
    \end{cases},
\end{equation*}
where the two species are self-promoted with mutual positive regulations, see Fig. 3 A. According to Corollary 1, the system will blow up within a finite time due to the positive Z-eigenvalues. Intuitively, self-promotion and positive regulation between the two species will result in an uncontrollable blow-up of the two populations very quickly, see Fig. 3 B. Based on Proposition 1, the blow-up time (i.e., singularity point) for the two species can be computed, which is given by $t=\min{\{1/(4 \alpha_1^2),1/(4 \alpha_2^2)\}}$ (where $\alpha_1$ and $\alpha_2$ depend on the initial conditions as described in Proposition 1). Although finite-time blow-up of populations is even worse than unbounded growth, it is actually possible to occur in modeling ecological systems \cite{parshad2013finite}. Hence, we can conclude that the two species cannot coexist under a confined ecological system.

\subsection{Population Model with Constant Control}
In this example, we consider a population model with pairwise and third-order interactions of the form:
\begin{equation}\label{eq:17}
    \frac{\dot{x}_i}{x_i}= r_i + \sum_{j=1}^n\sum_{k=1}^n \textsf{B}_{ijk}x_jx_k,
\end{equation}
where $x_i$ and $r_i$ represent the abundance and the intrinsic growth rate (i.e., reproduction and mortality rates) of species $s_i$, respectively, and $\textsf{B}$ is an interaction tensor that captures pairwise and third-order interactions among species. Detailed descriptions of third-order interactions can be found in \cite{bairey2016high}. If the second term in (\ref{eq:17}) is replaced by $\sum_{j=1}^n\textbf{B}_{ij}x_j$ (i.e., purely pairwise interactions), the population dynamics will become the classical generalized Lotka-Volterra model. For our purpose, we assume that the reproduction rate is equal to the mortality rate for each species, i.e., $r_i=0$. In fact, the population model described in the second example can also be modeled using (\ref{eq:17}). We further assume that each species has a supply rate (i.e., a constant input). Thus, we define our population model with pairwise and third-order interactions among three species as follows:
\begin{equation*}
\begin{cases}
    \dot{x}_1 &= -x_1^3-3x_1^2x_2-3x_1x_2^2 + 2\\
    \dot{x}_2 &= -x_2^3 + 2\\
    \dot{x}_3 &= -x_3^3-3x_1^2x_3-3x_1x_3^2-3x_2^2x_3\\ &-3x_2x_3^2-6x_1x_2x_3 + 2
\end{cases}.
\end{equation*}
In our model, each species is self-regulated (e.g., intraspecific competition) with a supply rate. The growth of species $s_1$ is inhibited by species $s_2$, and the growth of species $s_3$ is inhibited by species $s_1$, species $s_2$, and a combined effect from species $s_1$ and $s_2$.

The system of differential equations can be represented in the form of (\ref{eq:5}) such that $\textsf{A}\in\mathbb{R}^{3\times 3\times 3\times 3}$ is almost symmetric and $\textbf{b}=\begin{bmatrix}2 & 2 & 2\end{bmatrix}^\top$. Using Algorithm 1, we can confirm that the system can be transformed to an odeco HPDS with constant control, where the Frobenius norm error is about $1.27\times 10^{-15}$. The two factor matrices of \textsf{A} returned by Algorithm 1 are given by
\begin{equation*}
    \textbf{V}=\begin{bmatrix} 1 & 0 & 1\\1   &  1  &   1\\ 0  &   0  &   1\end{bmatrix} \text{ and } \textbf{V}^{(f)}=\begin{bmatrix}1  &   -1  &   0 \\ 0 &   1 &    0\\-1  &   0   & 1\end{bmatrix}
\end{equation*}
with weights $\lambda_1=\lambda_2=\lambda_3=-1$. For simplicity, we chose the standard basis to construct the odeco tensor, with which the transformation matrix is computed as
\begin{equation*}
    \textbf{P} =(\textbf{I} \textbf{V}^{-1})^\top=\begin{bmatrix} 1  &  -1  &   0\\ 0  &   1  &   0\\ -1  &   0  &   1\end{bmatrix}.
\end{equation*}
Suppose that $\textbf{x}(t)=\textbf{P}\textbf{y}(t)$ and $\textbf{y}(t)=\sum_{r=1}^3 c_r(t)\textbf{e}_r$ where $\textbf{e}_r$ are the three-dimensional standard basis vectors. Based on Proposition 4, the coefficient functions $c_r(t)$ can be solved implicitly by
\begin{equation*}
    \begin{cases}
    t &= \frac{g(2/3, 4/c_1^3)}{2c_1^2} - \frac{g(2/3, 4/\alpha_1^3)}{2\alpha_1^2}\\
    t &= \frac{g(2/3, 2/c_2^3)}{2c_2^2} - \frac{g(2/3, 2/\alpha_3)}{2 \alpha_2^2}\\
    t &= \frac{g(2/3, 6/c_3^3)}{2c_3^2} - \frac{g(2/3, 6/\alpha_3^3)}{2\alpha_3^2}
    \end{cases},
\end{equation*}
where $c_r(0)=\alpha_r$ for all $r$. According to the properties of Gauss hypergeometric functions, it can be shown that the three implicit functions have vertical asymptotes at $c_1=\sqrt[3]{4}$, $c_2=\sqrt[3]{2}$, and $c_3=\sqrt[3]{6}$, respectively, regardless of the initial conditions, see Fig. 4 A. This implies that the transformed system achieves global asymptotic stability at the equilibrium point. Thus, the equilibrium point of the original population dynamics, which is given by
\begin{equation*}
    \textbf{x}_e = \begin{bmatrix}0.3275 & 1.2599 & 0.2297\end{bmatrix}^\top,
\end{equation*}
is also globally asymptotically stable, see Fig. 4 B. Ecologically, we can conclude that the three species can coexist with supply rates.

\begin{figure}[t]
    \centering
    \includegraphics[scale=0.28]{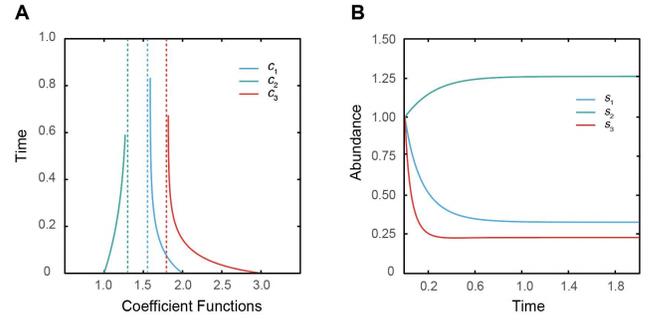}
    \caption{Population model with constant control example. (A) Plot of the three implicit coefficient functions. (B) Plot of the trajectories of the population dynamics. We used the initial condition $\textbf{x}_0=\begin{bmatrix}1 & 1 & 1\end{bmatrix}^\top$ in generating the two plots.}
    \label{fig:4}
\end{figure}

\section{Conclusion}\label{sec:7}
In this paper, we investigated the explicit solution and the stability properties of a continuous-time odeco HPDS. We derived an explicit solution formula using the Z-eigenvalues and the Z-eigenvectors from the orthogonal decomposition of the corresponding dynamic tensor. By utilizing the form of the explicit solution, the stability properties of the system could be formalized. In particular, the Z-eigenvalues can offer necessary and sufficient stability conditions. Furthermore, we explored the complete solution of an odeco HPDS with constant control. Finally, we provided criteria to determine if a general HPDS can be transformed to or approximated by an odeco HPDS with detailed algorithmic procedures. It will be worthwhile to investigate more strong results on approximating a general HPDS by an odeco HPDS. Future work also includes exploring stabilizability of odeco HPDS. For example, how to design a linear/polynomial control in order to shift the unstable Z-eigenvalues of an odeco HPDS to the left-half plane? Further, tensor algebra-based computation for Lyapunov equations and Lyapunov stability is important for future research.

\section*{Acknowledgements}
I would like to thank Dr. Anthony M. Bloch, Dr. Ram Vasudevan, Dr. Amit Surana, Dr. Andrea Aparicio, and Dr. Yang-Yu Liu for carefully reading the manuscript and for providing valuable comments. I would also like to thank the three referees for their constructive comments, which led to a significant improvement of the paper.

\bibliographystyle{IEEEtran}
\bibliography{reference}

\end{document}